\documentclass[]{amsart}

\usepackage{latexsym}
\usepackage{amssymb,amsmath,amsopn}
\usepackage[dvips]{graphicx}   
\usepackage{color,epsfig}      
\usepackage{url}
\usepackage{amsbsy}

\newtheorem{theorem}{Theorem}[]

\newtheorem{defi}{Definition}[]

\theoremstyle{definition}

\newtheorem{remark}[]{Remark}[]

\newtheorem{Pax}{Axiom}

\newtheorem{Sax}{Axiom}

\newtheorem{vyax}{Axiom}

\newtheorem*{axiomG}{Axiom G}

\newtheorem*{statement*}{Statement}

\begin{document}
\title[]{Gallucci's axiom revisited}
\author[\'A. G.Horv\'ath]{\'Akos G.Horv\'ath}
\address {Department of Geometry \\
Budapest University of Technology and Economics\\
H-1521 Budapest\\
Hungary}
\email{ghorvath@math.bme.hu}

\dedicatory{}

\date{Dec, 2017}

\begin{abstract}
In this paper we propose a well-justified synthetic approach of the projective space. We define the concepts of plane and space of incidence and also the Gallucci's axiom as an axiom to our classical projective space. To this purpose we prove from our space axioms, the theorems of Desargues, Pappus, the fundamental theorem of projectivities, and the fundamental theorem of central-axial collinearities, respectively. Our building up do not use any information on analytical projective geometry, as the concept of cross-ratio and the homogeneous coordinates of points.
\end{abstract}

\maketitle

\section{Introduction}

Very old and interesting question is the following: Is there a purely geometric reasoning of the fact that in the planes of a "three-dimensional projective space" the Pappus's theorem must be true? The classical investigations contains the following clear algebraic reasoning: A Desarguesian projective plane admits projective homogeneous coordinates of points in a division-ring. The constructed division-ring is a field if and only if the Pappus's theorem is true in it. Consequebtly, the Pappus' s theorem is necessary to define the cross-ratio of points and lines and also necessary to building up the standard tools of the so-called "classical projective geometry". Hence a classical projective plane can be defined as such a Desarguesian plane which the Pappus's theorem is valid. It is also known that every system of axioms to a higher-dimensional space implies the validity of the theorem of Desargues in its planes. However, in general, such a system of axioms says nothing about the validity of the theorem of Pappus; if the field of the coordinates is non-commutative then the Pappus's theorem fails to be valid. Since the commutativity is neede for the building up of the ususal tools of a classical projective space geometry the Pappus's axiom has to be assumed.

In this paper we prove that the Gallucci's axiom is a natural choice for an additional axiom to define the three-dimensional projective space. First of all it arises in its own right at the basic investigations of the transversals of skew lines. Its equivalence with the Pappus's theorem is known but most proofs uses either the concept of cross-ratio (see in \cite{kerekjarto}) or the involution theorem of Desargues (see in \cite{coxeter}). However, we could find in the book of Bachman, a nice figure suggesting a purely synthetic proof for this equivalence (see p.256 in \cite{bachman}). Bachmann investigated the nice theorem of Dandelin on "Hexagramme mistique". Dandelin generalised his theorem on conics to a hyperboloid of revolution, rather than a cone, relating Pascal's hexagon, Brianchon's hexagon and the hexagon formed by the generators of the hyperboloid. Dandelin's generalisation gives independent proofs of the theorems of Pascal and Brianchon.
This fact was proved in a metric geometry (using the axioms of incidence, order and congruence, respectively) by F. Schur in \cite{schur}.
Bachmann gives a reformulation of Dandelon's theorem which very similar to the Gallucci's statement, but without the simple formulation of this latter one. Coxeter raised in \cite{coxeter} that the statement of Gallucci (explicitely stated first in its form in \cite{gallucci}) can be used as an axiom of the projective space.

To realize this idea we give a short synthetic proof of this equivalence suggesting that Gallucci's statement. Additionally, we give an immediate proof of the equivalence of Gallucci's axiom and the fundamental theorem of projective mappings of  pencils of points. Finally, we also prove that the fundamental theorem of central-axial collinearities is also equivalent to the latter theorems. In this proof we do not use neither the fundamental theorem of projective geometry nor  any analytical tools from projective geometry.

\section{Incidence of the space-elements}

In the book of Veblen and Young \cite{veblen-young} we get a synthetic approach to define the alignment of an $n$-dimensional projective geometry which we call here the \emph{properties of incidences}. They used two undefined objects the point and the line, and define the concepts of plane, space, 4-space, etc. In this paper we consider these concepts also undefined consequently the simple system of assumptions used in the Veble-Young's we rewrite in a more didactic form. The mentioned assumptions of Veblen and Young in the three-case are the followings

\begin{vyax}
If $A$ and $B$ are distinct points, there is at least one line on both $A$ and $B$.
\end{vyax}
\begin{vyax}
If $A$ and $B$ are distinct points, there is not more than one line on both $A$ and $B$.
\end{vyax}
\begin{vyax}
If $A$, $B$, $C$ are points not all on the same line, and $D$ and $E$ ($D \ne E$) are points such that $B$, $C$, $D$ are on a line and $C$, $A$, $E$
are on a line, there is a point $F$ such that $A$, $B$, $F$ are on a line and also $D$, $E$, $F$ are on a line.
\end{vyax}
\begin{vyax}
There are at least three points on every line.
\end{vyax}
\begin{vyax}
There exists at least one line.
\end{vyax}
\begin{vyax}
All points are not on the same line.
\end{vyax}
\begin{vyax}
All points are not on the same plane.
\end{vyax}
\begin{vyax}
If $S_3$ is a three-space, every point is on $S_3$.
\end{vyax}

In our treatment the above "system of assumptions" will change to such statements which are theorems in the Veblen-Young system. We lost the requirement, of simplicity and independence but get a more didactic system of axioms leading to a faster building up of the three-dimensional geometry.

We call collinear (coplanar) some points if they are incident with the same line (with the same plane).  The incidence of points and lines may be considered in a so-called \emph{plane of incidence} which satisfies the following three axioms:

\begin{Pax}
Two points uniquely determine a line which is incident with them.
\end{Pax}
\begin{Pax}
Two lines uniquely determine a point which is incident with them.
\end{Pax}
\begin{Pax}
There exists four points, no three of which are collinear.
\end{Pax}

These concept can be found e.g. in the books \cite{reimann}, \cite{karteszi}.  For three-space we use the following set of axioms (see e.g. in \cite{gho}):

\begin{Sax}
Two points uniquely determine a line which is incident with them.
\end{Sax}
\begin{Sax}
Two planes uniquely determine a line which is incident with them. This is the line of intersection of the planes.
\end{Sax}
\begin{Sax}
Three non-collinear points uniquely determines a plane.
\end{Sax}
\begin{Sax}
Three planes which have no common line uniquely determine a point.
\end{Sax}
\begin{Sax}
If two points incident with a line are incident with a plane then all points of this line are also incident with this plane.
\end{Sax}
\begin{Sax}
If two planes incident with a line are incident with a point then all planes incident with this line are also incident with this point.
\end{Sax}
\begin{Sax}
In every planes there exists four points, no three of which are collinear.
\end{Sax}
\begin{Sax}
There exists five points, no four of which are coplanar.
\end{Sax}
We note, that from these axioms the usual set of assumptions (or axioms) of projective incidences can be proved. A similar coherent set of axioms can be found e.g. in the book of Ker\'ekj\'art\'o \cite{kerekjarto}, and the equivalence of the two system proved in \cite{gho}.

Clearly the planes of a space of incidence are plane of incidence, too, and it can be shown that the duality property of the projective space valid in every space of incidence.

\section{The role of the Desargues's theorem}

The plane of incidence is the pre-image of the concept of a classical projective plane. We can raise interesting problems in it, e.g. we can investigate that couple of statements which state the incidence of certain lines and point on the base of another and assumed incidences of the same points and lines. We refer to such a statement as a \emph{configuration theorem}. One of the most known example for configuration theorem is the theorem of Desargues. It says that
\begin{theorem}[Desargues]\label{thm:desargues}
The lines incident with the pair of points $A,A'$; $B,B'$ and $C,C'$ are go through the same point $S$ if and only if the points $X$, $Y$ and $Z$ incident with the pairs of lines $AB$, $A'B'$; $BC$, $B'C'$ and $AC$, $A'C'$ are incident with the same line $s$.
\end{theorem}

\begin{figure}[ht]
\includegraphics[scale=1]{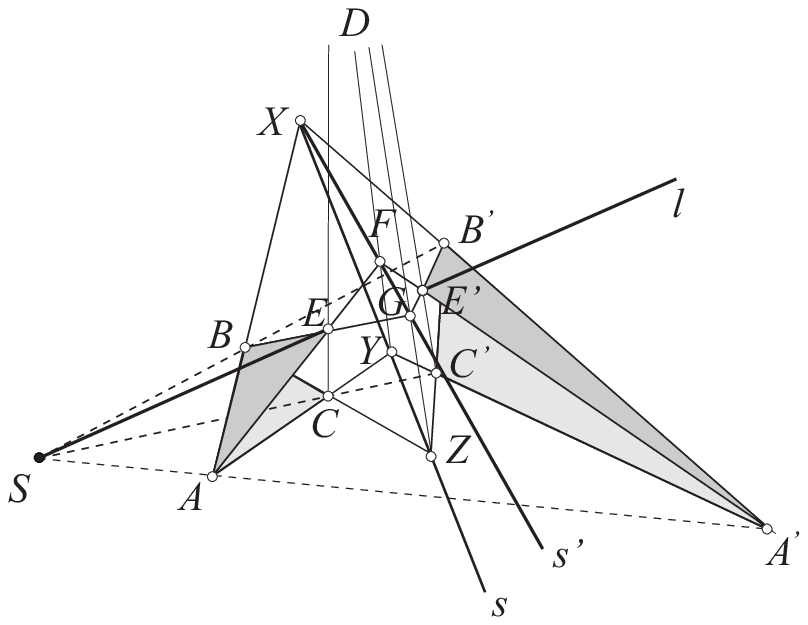}
\caption{The theorem of Desargues}
\label{fig:desargues}
\end{figure}
As a first question we can ask that without any further information on the plane of incidence we can prove or not the theorem of Desargues. The answer is known, in general this theorem does not true on a plane of incidence. (A nice counterexample has been given by Moulton in \cite{moulton},see also in \cite{gho}. For more information I propose the survey \cite{weibel} of Weibel.)

\begin{proof} Assume that the lines $AA'$, $BB'$ and $CC'$ of the plane $\alpha$ have a common point $S$. Since $\alpha$ is a plane of a space of incidence we can draw a line $l$ from $S$ which is not lying in $\alpha$. Consider that plane $\gamma $ which is defined by the intersecting lines $l$ and $CC'$ and consider a point $D$ of $\gamma $ which is not incident neither $l$ nor $CC'$. Let denote by $E=CD\cap l$ and $E'=C'D\cap l$, respectively. Both of the planes $ABE$ and $A'B'E'$ contains the point $X$ and distinct to each other hence they have a common line $s'$ which contains the point $X$. The lines $AE$ and $A'E'$ are lying in the respective planes $ABE$ and $A'B'E'$ and also in the plane $Al$ of the lines $AA'$ and $l$. Hence $AE$ is the common line of the planes $ABE$ and $Al$ and $A'E'$ is the common line of $A'B'E'$ and $Al$ then there exists an unique point $F$ which is the common point of the lines $AE$, $A'E'$ and $s'$. Similarly we can see that the lines $BE$, $B'E'$ and $s'$ are concurrent at their common point $G$. Consider now the planes $ACE$ and $A'C'E'$. The points $F$ and $D$ belong to these planes hence the line $DF$ is their common line. But $AC$ is the common line of $\alpha $ and $ACE$ and $A'C'$ is a common line of $\alpha $ and $A'C'E'$, hence the line $DF$, $AC$ and $A'C'$ meeting at their common point $Y$. Similarly, $DG$, $BC$ and $B'C'$ contains the point $Z$ implying that $X$, $Y$ and $Z$ are in the planes $\alpha $ and $Ds'$, respectively. Consequently these are on the common line of these planes, as we stated.

The converse statement follows from the proved direction if we consider the pairs of points $AB$, $A'B'$ and $YZ$ whose lines are intersect in the point $X$.
\end{proof}

\section{Transversals of lines in the space of incidence and the axiom of Gallucci}

A line is a \emph{transversal} of two other lines if it intersects both of them. In a space of incidence a pair of lines is called \emph{skew} if they are not lie in a common plane. The set of lines which are transversal of a pair of skew lines $a$ and $b$ form two pencils of planes with respective axes $a$ and $b$. Through every point $P$ of the space is lying neither on $a$ nor $b$, there is a plane through $a$ and an other through $b$. These planes intersect each other in an unique line $t$ through $P$. This line is called by the \emph{transversal of $a$ and $b$ from $P$}. If the space has infinitely many points, every three pairwise skew lines $a$, $b$, $c$ have infinitely many common transversals which are pairwise skew, respectively.

\begin{figure}[ht]
\includegraphics[scale=1]{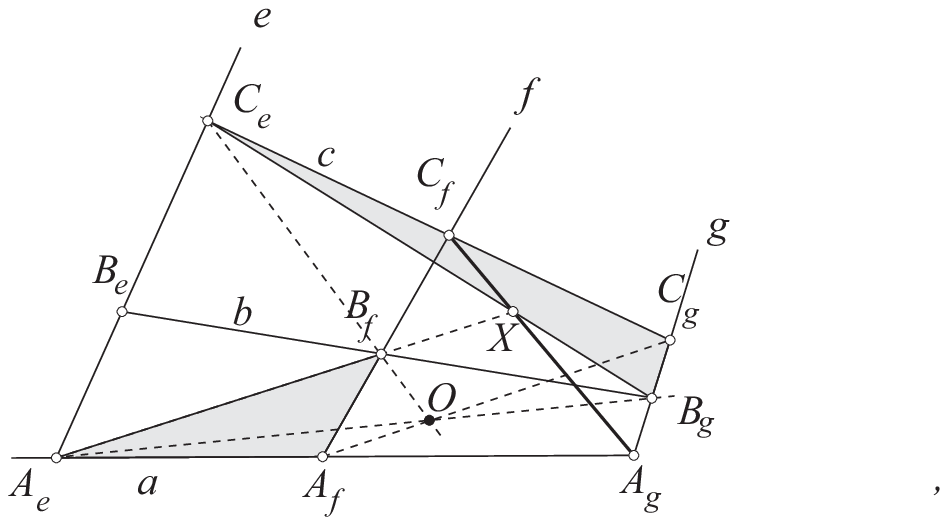}
\caption{The theorem of Desargues in space and the transversals of skew lines}
\label{fig:desarguesinthreespace}
\end{figure}

\begin{remark}
Note that the configuration of Desargues in the case of two non-coplanar triangles strongly connected to the existence of three common transversal $e,f,g$ of three pairwise skew lines $a,b,c$. Consider the figure Fig.\ref{fig:desarguesinthreespace} and its notation. Observe that the corresponding sides of the two triangles $(A_eA_fB_f)$ and $(B_gC_gC_e)$ are intersects to each other, if and only if $a:=(A_eA_f)$ and $g:=(B_gC_g)$ are intersecting, so the plane $\alpha=(a,g)$ is exist, $c:=(C_gC_f)$ and $f:=(A_f,B_f)$ are intersecting, hence $\beta:=(c,f)$ is exist and $A_eB_f$ and $B_gC_e$ are intersecting in a point $X$ equivalently, the plane $\gamma=(e,b)$ is also exist, where $e:=(A_eC_e)$ and $b:=(B_f,B_g)$. If either two from these planes are agree or the all incident with a line then the original two triangles are coplanar. From our assumption follows that we have three lines of intersection which meet in a point. These pairwise intersection lines are $(A_e,B_g)$, $(A_f,C_g)$ and $(B_f,C_e)$, respectively. which meet in a point $O$. Hence the two triangles are perspective with respect to $O$. On the other hand the existences of points $a\cap g$, $c\cap f$ and $e\cap b$ equivalent with their collinearity, because the respective pairs can meet on the line of intersection of the two triangle. As a conclusion we can state:
\emph{In a space of incidence two non-coplanar triangles are in perspective positions with respect to a point if and only if the corresponding sides intersect each other.}

From the figure Fig.\ref{fig:desarguesinthreespace} we can see immediately that every non-coplanar Desargues configuration (containing the ten points and ten lines of the theorem) can be associated to three pairwise skew lines $a,b,c$ and their transversals $e,f,g$. Two from the first group of lines (in the figure the lines $a$ and $c$) contain two skew edges of the triangles, respectively and the third line ($b=B_gB_f$) is the common line of the remaining vertices of the triangle; and two lines ($f$ and $g$) from the second group contain other two skew lines of the triangles, respectively and the third line $e$ is the common lines of the remaining vertices of the triangles.

\emph{Hence if we have two skew edges $a,c$ and their transversals $f,g$ and the line $b$ is a transversal of $f$ and $g$ moreover the line $e$ is a transversal of $a$ and $c$, the associated lines of Fig.\ref{fig:desarguesinthreespace} form an Desargues configuration if and only if the transversals $b$ and $e$ meet each other.}

Since all Desargues configuration is a projection of a non-coplanar Desargues one, we can associate any Desargues configuration a system of six lines, the first three lines transversal to the elements of the second group (and obviously vice versus).
\end{remark}

\begin{figure}[ht]
  \centering
    \includegraphics[scale=0.8]{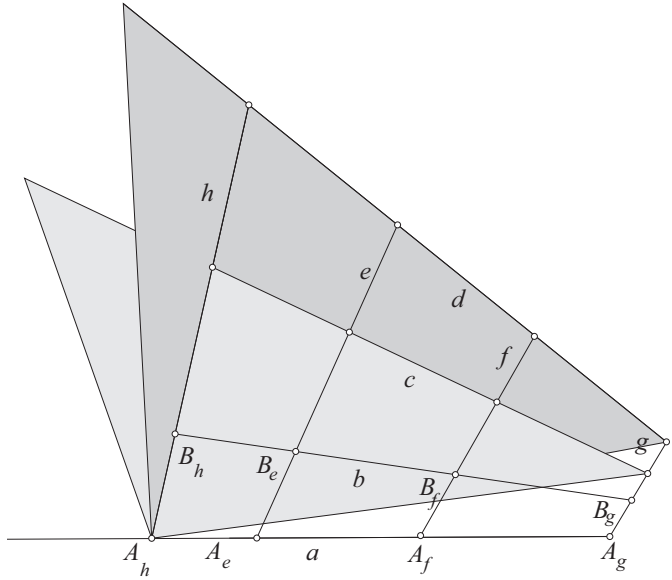}\\
  \caption{The Gallucci's axiom}\label{fig:gallucci}
\end{figure}

Consider now three pairwise skew lines $a,b,c$ and three common transversals $e,f$ and $g$, respectively. The plane $\gamma_e $ through the lines $c$ and $e$ intersects $a$ and $b$ in a corresponding pair of points $A_e$, $B_e$, respectively. (Consequently, $A_e,B_e$ and $C_e=c\cap e$ are collinear points.) Conversely if $A_e$ is any point of the line $a$ and we determine the plane $(c,A_e)$, the line $b$ intersects it a point $B_e$ determining the joining line $e:=(A_eB_e)$ of the plane $(c,A_e)$. This line $e$ intersects $c$ in a point $C_e$ which is collinear to $A_e$ and $B_e$. We can say for this situation that \emph{the point $A_e$ of $a$ is projected to the point $B_e$ of $b$ through the line $c$}. With this method we can project also $A_f$ to $B_f$ and $A_g$ to $B_g$ through $c$, respectively. If now $d$ is another common transversal of the lines $e,f,g$ then we can project the points of $a$ to the points of $b$ through the line $d$, too. Since $(c,e)\cap (d,e)=e$ the corresponding pairs of points at the two projections (through $c$ and through $d$, respectively) on the line $e$ are agree. The same holds in the cases of the lines $f$ and $g$, too. Consider now the planes $(c,A_h)$ and $(d,A_h)$ with a common point $A_h$ on $a$. \emph{If the line $b$ intersects these planes at the same point $B_h$} then the intersection of them is the line $(A_h,B_h)=:h$ intersects both of the lines $c$ and $d$, witnessing the so-called Gallucci's axiom:
\begin{axiomG}\label{ax:gallucci}
If three skew lines all meet three other skew lines, any transversal of the first set of three meets any transversal of the second set.
\end{axiomG}

If the Steiner-Pappus's theorem is valid in the space of incidence (consequently, we can speak about the concept of cross-ratio) the cross-ratio of the point $A_e,A_f,A_g$ and $A_h$ is equal to the cross-ratio of the corresponding planes through either $c$ or $d$ implying that the cross-ratios $(B_eB_fB_gB_h(c))$ and $(B_eB_fB_gB_h(d))$ are equal to each other. Since the cross-ratio uniquely determines the position of a point with respect to three fixed point we get that $B_h(c)=B_h(d)$ and the Gallucci's axiom is valid. (This is the case of the real projective space.)

\begin{remark}
We note that the statement of Gallucci's axiom is self-dual with respect to the natural duality of the space of incidence. In fact, the dual of skew lines are skew lines and the same  transversal line of two skew lines can be defined either the connecting line of two respective points of the lines or the intersection of two respective planes through the given lines.
So this statement in space may be similar fundamental property as the Desargues's theorem in the plane. We prove in the following that Gallucci's axiom implies all of the synthetic statement which needed to the building up of a projective space.
\end{remark}

\section{The role of the Pappus's axiom}

The Pappus's axiom is a very important statement of the real projective geometry. It is not valid in every projective plane  not even in every Desarguesian plane. Hessenberg proved in \cite{hessenberg} that in a plane of incident the Pappus's theorem implies the Desargues's theorem. This means that a Pappian plane is always Desarguesian, hence it has a coordinate field for the projective homogeneous coordinates of its points. In fact, this fields in the case of a Pappian plane is commutative.

\begin{theorem}[Pappus-Pascal]\label{thm:pp}
Assume that $A,B,C$ are three points of a line and $A',B',C'$ are three points of an other line. Then the points $C''=AB'\cap A'B$, $B''=AC'\cap A'C$ and $A''=BC'\cap B'C$ are collinear.
\end{theorem}

\begin{figure}[ht]
  \centering
    \includegraphics[scale=0.8]{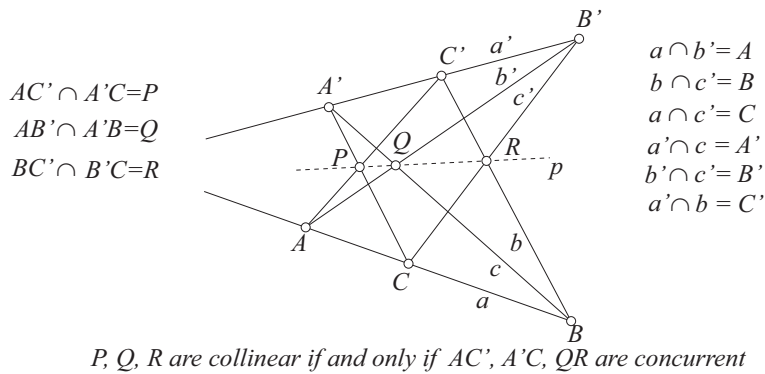}\\
  \caption{The equivalence of the Pappus's theorem and the theorem of Pappus-Brianchon}\label{fig:ppandpb}
\end{figure}

From the figure Fig. \ref{fig:ppandpb} we can check easily that the dual form of this theorem used first by Brianchon is equivalent to the above form. This latter says that

\begin{theorem}[Pappus-Brianchon]\label{thm:pb}
Assume that $a,b,c$ are three lines through a point and  $a',b',c'$ are three lines through another point. Then the lines $c''=(a\cap b',a'\cap b)$, $b''=(a\cap c',a'\cap c)$ and $a''=(b\cap c',b'\cap c)$ are concurrent.
\end{theorem}
In fact

\subsection{The Gallucci's axiom is equivalent to the Pappus's theorem in a space of incidence}

\begin{figure}[ht]
  \centering
    \includegraphics[scale=0.5]{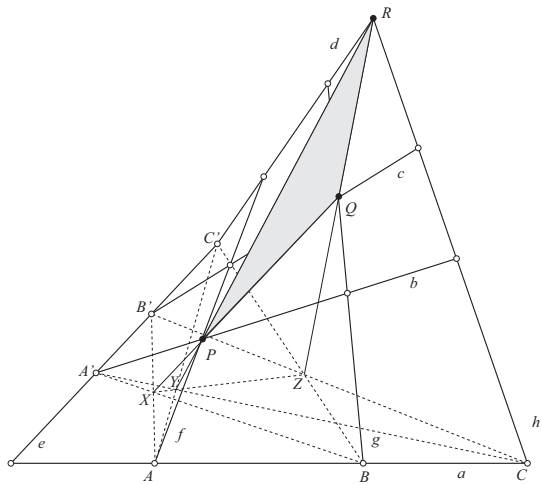}\\
  \caption{Gallucci's axiom equivalent to Pappus' axiom}\label{fig:gpp}
\end{figure}

On the figure Fig.\ref{fig:gpp} we can see the connection between the Gallucci's axiom and Pappus's theorem. 

\noindent{\bf Pappus' axiom implies Gallucci's axiom:}
On the figure Fig.\ref{fig:gpp} we can see the connection between the Gallucci's axiom and Pappus' axiom. Assume first that the Pappus' axiom is valid in our planes and consider three pairwise skew lines $a$,$b$ and $c$ with three transversals $e$, $f$, and $g$, respectively. Assume that $d$ is a transversal of $e,f,g$ and $h$ is a transversal of $a,b,c$, respectively. Then we have to prove that $d$ intersects $h$ in a point.
By the notation of Fig.\ref{fig:gpp} in the plane of the intersecting lines $a$ and $e$ we get a configuration of nine points $A,B,C$; $A',B',C'$, and $X=AB'\cap A'B$,$Y=AC'\cap A'C$,$Z=BC'\cap B'C$. If the axiom of Pappus is valid then the points $X,Y,Z$ are collinear. Let $P=b\cap f$ and $Q=c\cap g$ and consider the plane $\pi$ of the points $X,Y,Z$ and $P$. Since the plane of $b$ and $g$ contains the points $P,A',B$ it contains the point $X$ and also the line $PX$. Similarly the plane of $c$ and $f$ contains the points $P,A,B'$ and hence the point $X$ and the line $PX$, respectively. From this we get that $P,X$ and $Q$ are collinear points. Consider the intersection $R$ of the coplanar lines $YP$ and $ZQ$. Since the point $Z$ is on the line $CB'$, and $Q$ is on the line $c$, the line $ZQ$ is on the plane of the lines $c$ and $B'C$ which is the same as the plane of $c$ and $h$. Similarly, the line $YP$ is on the plane of the points $A',C,P$ which is the plane of the lines $b$ and $h$. Then the intersection point $R$ of the lines $ZQ$ and $YP$ has to lie on the intersection of the planes $(c,h)$ and $(b,h)$ which is the line $h$. Similarly, the line $ZQ$ is in the plane $(d,g)$; and the line $YP$ is in the plane $(d,f)$, showing that the point $R$ is also in the line $d$. Hence $d$ and $h$ intersect each other in a point as the Gallucci's axiom states.

\noindent{\bf Gallucci's axiom implies Pappus' theorem:}
To prove that the Gallucci's axiom implies the Pappus' theorem consider a Pappian configuration of the nine points $A,B,C,A',B',C',X,Y,Z$ of Fig.\ref{fig:gpp} in the plane of the lines $a$ and $e$. We would like to prove that the points $X,Y$ and $Z$ are collinear points. To this purpose consider a point $P$ of the space of incidence which is not lying on the plane $(a,e)$ and define the lines $b:=A'P$ and $f:=AP$, respectively. Let $Q$ be an arbitrary point of the line $XP$ distinct from $X$ and $P$, and define the lines $c:=B'Q$ and $g:=BQ$, respectively. Finally, let $d$ be the unique transversal from $C'$ to the skew lines $f$ and $g$ and let $h$ be the unique transversal from $C$ to the skew lines $b$ and $c$, respectively. If the Gallucci's axiom is true then $d$ and $h$ intersect each other in a point $R$. By the definition of the point $Y$ it is on the plane $(b,h)$ and it is also on the plane $(f,d)$ hence the line $YP$ is the line of intersection of these two planes. So the common point $R$ of $d$ and $h$ should lie on $YP$. Similarly, the points $Z$ and $Q$ lie on the planes $(c,h)$ and $(g,d)$ respectively, so $Z$, $Q$ and $R$ are collinear points, too. Hence the points $X,Y,Z$ are on the intersection of the planes $(a,e)$ and $(P,Q,R)$ implying that they are collinear points, as we stated.

\subsection{Gallucci's axiom is equivalent to the fundamental theorem of projectivities}

Given two lines $a$ and $b$ in a plane and a point $P$ of that plane on neither line, the bijective mapping between the points of the range of $a$ and $b$  determined by the lines of the pencil on $P$ is called a \emph{central perspectivity} (or more precisely, a central perspectivity with center $P$). A finite product of perspectivities is called by \emph{projectivity} of the space of incidence. The fundamental theorem of projectivities says that a projectivity is determined by the image of three points of a projective range (the set of points on a line). It is known that this statement is equivalent to the Pappus' theorem (see e.g. \cite{veblen-young}). Using cross-ratio the first direction of the equivalence is trivial but to prove the second one is a more complicated task. In this paper we avoid the concept of cross-ratio and prove this equivalence immediately.

To this purpose observe first that Gallucci's axiom says that the perspectivity of a point row to another point row through a line is independent of the choice of the axes of perspectivity. More precisely if $a,b,c$ are three pairwise skew lines and $e,f,g,h...$ common transversals of these, then the correspondence $a\cap e\mapsto b\cap e$, $a\cap f\mapsto b\cap f$, $a\cap g\mapsto b\cap g$, $a\cap h\mapsto b\cap h$, and so one... defines a mapping $\Phi(c)$ on the projective range of $a$ to the projective range of $b$ through the line $c$. This mapping is determined by the line $c$. Let now $d$ be any common transversal of the lines $e,f,g$ distinct from $a,b,c$. Similarly, from the points of $a$ to the points of $b$ there is a similar mapping $\Phi(d)$ through $d$ which agrees with $\Phi(c)$ in the three points $a\cap e$, $a\cap f$ and $a\cap g$. Let now $h$ be a common transversal of $a,b,c$ distinct from $e,f,g$, respectively. Gallucci's axiom is true if and only if the lines $d$ and $h$ intersect each other in a point, implying that the equality $\Phi(c)(a\cap h)=\Phi(d)(a\cap h)$ holds for every point $a\cap h$ of $a$. Hence any common transversal $d$ of $e,f,g$ can be the axis of the perspectivity under the investigation.

On the other hand the perspectivity of the point row $a$ to the point row $b$ through the line $c$ is the product of two central perspectivities. In fact, if we project from a point $C$ of the line $c$ the points of $b$ to a plane $\alpha$ through the line $a$, we get a projective range of a line $b'$ denoted by $(b\cap e)'$ $(b\cap f)'$, etc..., respectively. These points are perspective images of the respective points $a\cap e$, $a\cap f$, etc... of the line $a$ through the point of $\alpha\cap c=O(c)$. Hence the mapping $\Phi(C)$ is the product of the perspectivity $a\cap e\mapsto (b\cap e)'$ through $O(c)$ and the perspectivity $(b\cap e)'\mapsto b\cap e$ through $C$. Similarly, the mapping $\Phi(d)$ is also the product of two perspectivity with centers $O(d)$ and $D$.

Hence the two perspectivities through the respective lines $c$ and $d$ coincide if and only if the two products of perspectivities through the pairs of centers $O(c),C$ and $O(d),D$ coincide. By definition these two products agree in the points $a\cap e$, $a\cap f$ and $a\cap g$ because of $c$ and $d$ are transversals of the lines $e,f,g$. Consequently the two perspectivities through the lines $c$ and $d$ coincide if the fundamental theorem of projectivities of point row is valid in our space of incidence.

To prove that the Gallucci's theorem implies the fundamental theorem of projectivities first we prove that any projectivity of the skew lines $a$ and $b$ can be considered as a perspectivity of the point row $a$ to the point row $b$ through the line $c$.

First observe, that if the lines $a$, $a'$ and $a''$ meet in a common point, then the product of the two perspectivities $a\rightarrow a'$ $A\mapsto A'$ with center $O$ and $a'\rightarrow a''$ $A'\mapsto A''$ with center $O'$ can be simplified into one perspectivity $a\rightarrow a''$ with such a center $O^\star$ which lies on the line $OO'$. To prove this, we have to use the Desargues' theorem implying that the lines $AA''$, $BB''$, $CC''$ ... go through in the same point $O^\star$ of the line $OO'$.

As a second observation, we note that if the lines $a$, $a'$ and $a''$ have no common point, but $a\cap a'\ne \emptyset$ and $a'\cap a''\ne \emptyset$, then we can change the point row $a'$ into any other point row $b$, which go through the point $a\cap a'$, and does not contain the center $O'$. More precisely, we can change the second perspectivity $a'\overset{O'}{\rightarrow} a''$ to the composition $a'\overset{O'}{\rightarrow} b \overset{O'}{\rightarrow} a''$. Then the investigated mapping $a\overset{O}{\rightarrow}a'\overset{O'}{\rightarrow}a''$ can be considered as the new product $a\overset{O}{\rightarrow}a'\overset{O'}{\rightarrow}b\overset{O'}{\rightarrow}a''$. By our first observation it can be simplified to the form
$a\overset{O^\star}{\rightarrow}b\overset{O'}{\rightarrow}a''$.

Consider now a finite sequence of perspectivities $a^{i}\overset{O^{i}}{\rightarrow}a^{i+1}$ where $i$ runs from $1$ to $n-1$, then we can simplify the representation of the given projectivity as follows. We can assume that there is no three consecutive perspectivities with concurrent axes. If we have at least three perspectivities in the sequence that we choose the line $b$ of the previous observation to the line is determined by the points $a^1\cap a^2$ and $a^3\cap a^4$. $b$ evidently avoid the center $O^2$ and we can use the second observation to the first pair of the projectivities with this $b$. Then we get that the original product $a^1\overset{O^1}{\rightarrow}a^2\overset{O^2}{\rightarrow}a^3\overset{O^3}{\rightarrow}a^4$ can be changed to the following one:
$a^1\overset{O^\star}{\rightarrow}b\overset{O^2}{\rightarrow}a^3\overset{O^3}{\rightarrow}a^4$ where the lines $b$, $a^3$ and $a^4$ are concurrent. This means that using the first observation we can decrease the number of perspectivities. From an inductive argument this simplification leads to a product of two (suitable) perspectivity. If $a^1$ and $a^n$ are skew, that the simplified chain of two perspectivities is of the form $a^1\overset{O}{\rightarrow}b^{n-1}\overset{O^n}{\rightarrow}a^n$. Clearly, all lines $A^1A^n$, $B^1B^n$, $C^1C^n$,... intersect the line $OO^n$. (Observe that e.g. the line $A^1X^{n-1}$ contains the center $O$ and the line $A^nX^{n-1}$ contains the center $O^n$, respectively, implying that $A^1A^{n}$ and $OO^n$ are coplanar.) Hence the projectivity is the perspectivity of $a^1$ to $a^n$ through the line $OO^n$. Now if we have two projectivities from $a=a_1$ to $b=a_n$ which send the point $A^1$ to $A^n$, the point $B^1$ to $B^n$ and the point $C^1$ to $C^n$, respectively, then they can be interpreted as two perspectivities from $a$ to $b$ through the respective lines $c$ and $d$. These two maps are agree at three points of the line $a$. If Gallucci's theorem true then the image of any further point of $a$ is the same to the two mappings implying that the two projectivities coincide.

The last note that if the original lines having in a common plane that we can compose the respective projectivities by the same projectivity sending the second line to another one which is skew with respect to the first one. Using again the above argument we get that the original projectivities also coincide and the fundamental theorem of projectivities are true.

\subsection{The fundamental theorem of central-axial collinearities}

A \emph{collinearity} of the projective plane to itself is a bijective mapping of its points with the property that the image of a line is also a line. The collinearity is \emph{central} (\emph{axial}) if there is a point (line) of the plane with the property, that the image of a line (point) through this point (line) is equal to itself. It can be proved, that every central collinearity is an axial one and vice versa. From the theorem of Desargues it can be proved that a central-axial collinearity (shortly c-a collinearity) is determined by it centre, its axis, and a pair of points from which the second is the image of the first by this mapping. Obviously, the composition of central-axial collinearities are not central-axial collinearities. We consider now a finite product of central-axial collinearities, which we call a projective mapping of the plane to itself.  Thus arise the problem that four general points and their images determine or not a projective mapping of the plane. The fundamental theorem of c-a collinearities says that the answer is affirmative. 

We prove that {\bf the fundamental theorem of projectivities implies the fundamental theorem of c-a collinearities}.
First we show that if the pairs of points $A,A'$, $B,B'$, $C,C'$ and $D,D'$ are in general position then there is at least one projectivity $\phi$ of the plane with the property: $A'=\phi(A)$, $B'=\phi(B)$, $C'=\phi(C)$ and $D'=\phi(D)$. This projectivity of the plane is the product of at most four c-a collinearities.

\begin{figure}[ht]
  \centering
    \includegraphics[scale=1]{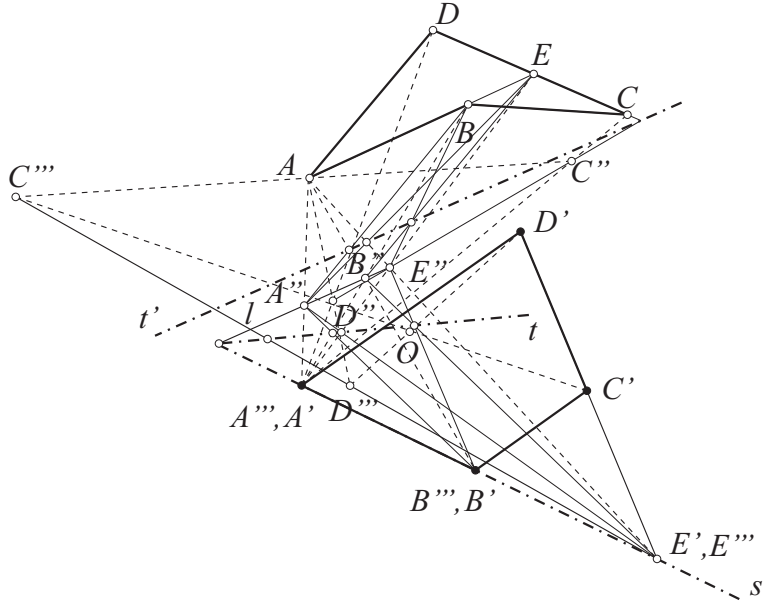}\\
  \caption{Composition of c-a collinearities can send given four points two another given four points}\label{fig:cacoll}
\end{figure}

Assume that the line $AB$ is not the line $A'B'$. In this case the points $A,B,E:=AB\cap CD$ and $A',B',E':=A'B'\cap C'D'$ form two triplets of two distinct lines as we can see in Fig.\ref{fig:cacoll}. We give two c-a collinearities which product sends $A$ to $A'$; $B$ to $B'$ and $E$ to $E'$ respectively. This product sends $C$ to $C'''$ and $D$ to $D'''$, respectively. Since $C'''D'''\cap C'D'= E'$ the third collinearity $\phi_3$ with center $O:=C'C'''\cap D'D'''$ and axis $A'B'$ with a pair of points $C''',C'$ fixes $A',B'$ and sends $D'''$ to $D'$, implying the required result. We have to give now the first two c-a collinearities, respectively. First consider the line $l$ of the two points $B''=AB'\cap A'B$ and $E''=AE'\cap E'A$. The intersection of $l$ with the line $AA'$ denote by $A''$. Observe that the points $A,B,E$ are in perspective positions with the points $A'',B'',E''$ from the point $A'$. Using Desargues' theorem to the triplets $A,B'',E$ and $A'',B,E''$ we get that the three lines $AB$, $l$ and $t':=(AB''\cap A''B,BE''\cap B''E)$ go through the same point. Similarly, the same true for the lines $AB$, $l$ and $(AB''\cap A''B,AE''\cap A''E)$ implying that third line in the two triplets is the same. Hence the c-a collinearity $\phi_1$ defined by the center $A'$, axis $t'$ and corresponding pair of points $A, A''$ sends $A,B,E$ to the points $A'',B'',E''$, respectively. Set $\phi_1(C):=C''$ and $\phi_1(D):=D''$. 

Define the second c-a collinearity $\phi_2$ on an analogous way, originated from the points $A'',B''$ and $E''$ and getting as the respective image points $A'$, $B'$ and $C'$. The center of $\phi_2$ is $A$ the axis of it is $t:=(A''B'\cap A'B'',B''E'\cap B'E'')$ and a pair of its points is $A'',A'$. If we denote by $C'''$ and $D'''$ the image at $\phi_2$ of the points $C''$ and $D''$ we get that situation from which we defined $\phi_3$.

If finally the line $AB$ is the same as the line $A'B'$ we apply a suitable c-a collinearity $\phi_0$ which sends the line $AB$ in another line $A_0B_0$ and for this situation we use the above arguments getting the third mappings $\phi_1$, $\phi_2$ and $\phi_3$. Now the product of the four c-a collinearity gives the required result.

Assume now that there are two projectivities (say $\phi$ and $\psi$) which sends the points $A,B,C,D$ into the points $A',B',C',D'$, respectively. Then the projectivity $\phi\circ \psi^{-1}$ gives a projectivity which fixes the points $A,B,C,D$ and distinct from the identity. Since the point $E=AB\cap CD$ are on two invariant lines, it is also a fix point. Since the fundamental theorem of projectivities of the point rows are true all the points of these two lines are fixed. Let $P$ be any point of the plane having outside from the lines $AB$ and $CD$. Clearly we have two lines with common point $P$ which intersects $AB$ and $CD$ implying that the point $P$ is also a fixed point of the projectivity. Hence the mapping is the identity, contrary to our indirect hypothesis.

Conversely we prove on the Desarguesian plane that {\bf the fundamental theorem of c-a collinearities implies the fundamental theorem of projectivities}. As we saw in the first part of this section, a perspectivity between two lines $a$ and  $b$ can be considered as the impact of a $\phi_{a,b}$ c-a collinearity on $a$, where the center of the perspectivity is the center of the c-a collinearity. Hence a finite product of perspectivities can be considered as the impact of a finite product of c-a collinearities on the given line $a$. If $A$ and $B$ two points of $a$ and $C$, $D$ two such points for which the four point lie in general position, then there exists the point $E=AB\cap CD$ and $A$, $B$ and $E$ are distinct points. Consider now two projectivities $\phi $ and $\psi$ on the point row $a$ to the point row $a'$ which are sending $A$ to $A'$, $B$ to $B'$ and $E$ to $E'$. The corresponding products $\phi^\star$ and $\psi^\star$ of c-a collinearities are sending $C,D$ to $C'''$, $D'''$ and $C,D$ to $C''''$, $D''''$, respectively. Define the c-a collinearity $\eta^\star$ with the axis $a'$, center $O=C'''C''''\cap D'''D''''$ and the equalities $\eta(C''')=C''''$ and $\eta(D''')=D''''$. (Of course,
if the line $C'''D'''$ distinct from  the line $C''''D''''$ the center is the intersection of the lines $C'''C''''$ and $D'''D''''$, otherwise it is a point of this common line can be determined easily, too.) We have $\eta^\star\circ \phi^\star=\psi^\star$ implying that the impact of the two projectivities are agree on the line $a$, too. Since $a'$ is a point-wise fixed line at $\eta^\star$, we proved that $\phi$ and $\psi$ are agree all of the points of $a$, too. Hence the fundamental theorem of projectivities of point rows are true on the plane, as we stated.

\section{Conclusion}

By the content of the previous sections we suggest the following definition of the three-dimensional projective space:

\begin{defi}
A space of incidence is projective space of dimension three, if the Gallucci's axiom holds in it. More precisely, a system of points, lines and planes form a projective space if the axioms ${\bf S1}-{\bf S8},{\bf G}$ are hold.
\end{defi}

As we saw in a projective space of dimension three the Pappus-Pascal theorem is valid implying that the division ring of the coordinates of its points is a commutative field. Clearly the cross-ratio of four points of a line can be defined correctly, and we can use the analytic methods of the projective geometry, too.

\end{document}